\documentclass[11pt, leqno,twoside]{article}
\usepackage{amssymb}
\usepackage{amsmath}
\usepackage{amsthm}
\usepackage{amsfonts}
\usepackage{color}
\usepackage{textcomp}
\usepackage{graphicx}

\allowdisplaybreaks

\begin{document}

\title{A note on Liouville theorem for  stationary flows of shear thickening fluids in the plane 
\thanks
{E-mail addresses: \   guozhang@jyu.fi}
}
\author{{ $\text{ Guo Zhang}$} \\
{\small Department of Mathematics and Statistics, P.O. Box 35 (MaD), FI-40014}\\
{\small University of Jyv\"{a}skyl\"{a}, Finland }}

\date{}
\maketitle

\begin{abstract}
In this paper we consider the entire weak solutions of the equations for stationary flows of shear thickening fluids in the plane and prove
Liouville theorem under the global boundedness condition of velocity fields.

 \end{abstract}
{\bf MR Subject Classification:}\ \ 76 D 05, 76 D 07, 76 M 30, 35 Q 30.
\bigbreak

\begin{center}
\section*{ \S1.\ Introduction }
\end{center}

In this note, we study entire weak solutions $u:\mathbb{R}^2 \rightarrow \mathbb{R}^2,\  \pi:\mathbb{R}^2 \rightarrow \mathbb{R}$ of the following system
\[
\left\{
\renewcommand{\arraystretch}{1.25}
\begin{array}{lll}
- \text{div}[T(\varepsilon(u))]+u^k\partial_{k}u+D\pi=0,\\
 \text{div}\ u=0\ \ \ \ \text{in}\ \  \mathbb{R}^2
\end{array}
\right.
\eqno(1.1)
\]
and show that bounded solutions are constants.

The above system describes the stationary flow of an incompressible generalized Newtonian fluid, where $u$ is the velocity field, $\pi$
the pressure function,  $u^k\partial_{k}u$  the convective term, and $ T$ represents the stress deviator tensor. As usual $\varepsilon(u)$ stands
 for the symmetric part of the differential matrix $Du$ of $u,$ i.e.

\[
\renewcommand{\arraystretch}{1.25}
\begin{array}{lll}
\varepsilon(u)=\frac{1}{2}(Du+(Du)^{T})=\frac{1}{2}(\partial_iu^k+\partial_ku^i)_{1\leq i,k \leq 2}.
\end{array}
\]

we assume that the stress tensor $T$ is the gradient of a potential $H:S^{2\times 2}\rightarrow \mathbb{R}$ defined on the space $S^{2\times 2}$ of
 all symmetric $(2\times 2)$ matrices of the following form
\[
\renewcommand{\arraystretch}{1.25}
\begin{array}{lll}
H(\varepsilon)=h(|\varepsilon|),
\end{array}
\eqno(1.2)
\]
where $h$ is a nonnegative function of class $C^2$. Thus
\[
\renewcommand{\arraystretch}{1.25}
\begin{array}{lll}
T(\varepsilon)=DH(|\varepsilon|)=\mu(|\varepsilon|)\varepsilon, \ \ \ \mu(t)=\frac{h'(t)}{t}.
\end{array}
\eqno(1.3)
\]

If $\mu$ is a constant, $\mu=\nu$, that is, if $h(t)=\frac{\nu}{2}t^2$, then system $(1.1)$ reduces to the stationary Navier-Stokes
equations for incompressible Newtonian fluids with viscosity coefficient $\nu$.

If $\mu$ is not a constant, then it means that the viscosity
coefficient depends on $\varepsilon$, and system $(1.1)$ describes
the motion of continuous media of generalized Newtonian fluids. For
the physical background and mathematical theory of generalized
Newtonian fluids, we refer to Ladyzhenskaya \cite{La}, Galdi
\cite{Ga1,Ga2}, Malek, Necas, Rokyta and Ruzicka \cite{MNRR}, and
Fuchs and Seregin \cite{FS}.

In the whole paper, we will concentrate on some special types of
shear thickening fluids. We assume that  the potential $h$
satisfies the following conditions:
\[
\renewcommand{\arraystretch}{1.25}
\begin{array}{lll}
h\  \text{is strictly increasing and convex}\\
\text{together with $ h''(0) > 0$ and $\lim\limits_{t\rightarrow 0}\frac{h(t)}{t}=0$};
\end{array}
\eqno(A1)
\]
\[
\renewcommand{\arraystretch}{1.25}
\begin{array}{lll}
\text{(doubling property) there exists a constant $a\geq 1$} \\
\text{ such that} \  h(2t)\leq a h(t)\  \text{for all} \  t \geq 0;
\end{array}
\eqno(A2)
\]
\[
\renewcommand{\arraystretch}{1.25}
\begin{array}{lll}
\text{we have}\  \frac{h'(t)}{t} \leq h''(t)\  \text{for any}\  t \geq 0.
\end{array}
\eqno(A3)
\]

The study of Liouville type theorems goes back to the work of Gilbarg and Weinberger \cite{GW}. They showed that entire solutions $u$ of
stationary Navier-Stokes equations in the plane are constants, provided that $\int_{\mathbb{R}^2}|Du|^2 dx  < \infty.$ For the unstationary
Navier-Stokes equations in $2D$, recently, Koch, Nadirashvili, Seregin and Sverak  \cite{KNSS} showed
that  $u(x,t)=b(t)\  \text{on}\  \mathbb{R}^2\times (-\infty,0)$ provided the solutions are bounded. Clearly, this
result implies the Liouville theorem for stationary Navier-Stokes equations, that is, bounded solutions to stationary Navier-Stokes
equations are constants.

In the case of shear thickening fluids, for $h$ satisfying $(A1)-(A3)$, Fuchs \cite{Fu1} very recently proved the following Liouville theorem.

 \textbf{Theorem }[Theorem $1.2$, \cite{Fu1}]\ \ Let \ $u\in C^2(\mathbb{R}^2,\mathbb{R}^2)$, $\pi\in C^1(\mathbb{R}^2,\mathbb{R})$ be the
 solutions to $(1.1)$. Suppose that $u$ is bounded in $\mathbb{R}^2$ and satisfies that
\[
\renewcommand{\arraystretch}{1.25}
\begin{array}{lll}
\sup\limits_{\mathbb{R}^2-B_R(0)}|u-u_{\infty}|\rightarrow 0\
\end{array}
\eqno(\ast)
\]
as $R \rightarrow \infty$, for a vector $u_{\infty}\in \mathbb{R}^2. $ Then $u$ is a constant vector.

In \cite{Fu1}, Fuchs conjectured that one can remove the assumption $(\ast)$ on $u$ at infinity in the above theorem and show
that any bounded solution $u$
must be a constant vector. In this note, we will give a positive answer to this conjecture.

 \textbf{Theorem 1}\ \ Suppose $u\in C^1(\mathbb{R}^2,\mathbb{R}^2) \bigcap L^{\infty}(\mathbb{R}^2,\mathbb{R}^2)$ be an entire
 weak solution to $(1.1),$ i.e.
\[
\renewcommand{\arraystretch}{1.25}
\begin{array}{lll}
 \displaystyle \int_{\mathbb{R}^2}T(\varepsilon(u))\colon \varepsilon(\varphi) dx - \displaystyle \int_{\mathbb{R}^2}u^k u^i \partial_k \varphi^i dx = 0
\end{array}
\eqno(1.4)
\]
for all $\varphi \in C_0^\infty(\mathbb{R}^2,\mathbb{R}^2), \text{div}\varphi=0.$ Then $u$ is a constant vector.

We first comment on the  regularity assumption on the  solutions.  As we known, for general $h$ satisfying $(A1)-(A3)$, $C^{1,\alpha}$ regularity
of solutions of $(1.1)$ is an open problem, but in some special case, such as $h(t)=t^2(1+t)^m, m\geq 0,$  the solution $u$ belongs to the
space $C^{1,\alpha}$ \cite{BFZ}, so the assumption in Theorem $1$ is reasonable.

Second, we comment on the proof of Theorem $1$. Our proof follows the same line as that of Theorem $1.2$ in \cite{Fu1}. We first need
an energy estimate for the first order derivatives, as stated in Lemma $3.1$. This was proved in \cite{Fu1}. Then we need the energy estimate
for the second order derivatives, see formula $(3.10)$, from which follows that
\[
\renewcommand{\arraystretch}{1.25}
\begin{array}{lll}
\displaystyle\int_{\mathbb{R}^2} D^2H(\varepsilon(u))(\varepsilon(\partial_ku), \varepsilon(\partial_ku))dx < \infty.
\end{array}
\eqno(\star)
\]
Actually, $(\star)$ was proved in \cite{Fu1}, via a Caccioppoli-type
inequality, see $(5.10)$ in \cite{Fu1}. Our Caccioppoli-type
inequality is different from  $(5.10)$ in \cite{Fu1}. For the proof
of $(3.10)$,  we use a different approach than that in  \cite{Fu1}
to estimate the item involving the pressure $\pi$. The treatment is
standard. The essential point of our proof is based on the delicate
analysis of Caccioppoli-type inequality $(3.10)$. We estimate one
crucial item involving the first order derivatives by integration by
parts. This idea was  also used in \cite{FZ}.

 Our notations is standard. Throughout this paper,  the convention of summation with respect to indices repeated twice is
 used. All constants are denoted by the symbol $ C,$ and  $ C $ may change from line to line.  Whenever it is necessary we
 will indicate the dependence of
$C$ on parameters. As usual $Q_R(x_0)$ denotes the open square with center $x_0$ and side length $2R,$  and
symbols $\colon,$  $\cdotp$ will be used
 for the scalar product of matrices and vectors respectively. $|\cdotp|$ denotes the associated Euclidean norms.

Our paper is organized as follows: in Section $2$ we present some
auxiliary results, and in Section $3$ we give the proof of Theorem $1$.

\bigbreak
\begin{center}
\textbf{\Large {\S 2.\ Auxiliary Results}}
\end{center}

As shown in \cite{Fu1}, the following properties of functions $h$ follows from $(A1)-(A3)$.

$(i)$\ \ $\mu(t)=\frac{h'(t)}{t}$ is an increasing function.

$(ii)$ \ \ We have $h(0)=h'(0)=0$ and
\[
\renewcommand{\arraystretch}{1.25}
\begin{array}{lll}
h(t) \geq \frac{1}{2} h''(0)t^2.
\end{array}
\eqno(2.1)
\]

Moreover,
\[
\renewcommand{\arraystretch}{1.25}
\begin{array}{lll}
\frac{ h'(t)}{t} \geq \lim\limits_{s\rightarrow 0}\frac{h'(s)}{s}=h''(0)>0.
\end{array}
\eqno(2.2)
\]

$(iii)$\ \ The function $h$ satisfies the balancing condition, i.e.
\[
\renewcommand{\arraystretch}{1.25}
\begin{array}{lll}
\frac{1}{a}h'(t)t\leq h(t) \leq th'(t),\ \ \ \ t\geq 0.
\end{array}
\eqno(2.3)
\]

$(iv)$\ \ For an exponent $m \geq 2$ and a constant $C \geq 0$ it holds
\[
\renewcommand{\arraystretch}{1.25}
\begin{array}{lll}
h(t)\leq C (1+t^m),\ \ \ \ h'(t)\leq C(1+t^m),\ \ \ \ t \geq 0.
\end{array}
\eqno(2.4)
\]

From the assumptions on $h,$ we know the system satisfies the following elliptic condition, $\forall \varepsilon, \sigma \in S^2,$
\[
\renewcommand{\arraystretch}{1.25}
\begin{array}{lll}
\frac{h'(|\varepsilon|)}{|\varepsilon|}|\sigma|^2 \leq D^2H(\varepsilon)(\sigma,\sigma)\leq h''(|\varepsilon|)|\sigma|^2,
\end{array}
\eqno(2.5)
\]

from which, together with $(2.2)$, follows that

\[
\renewcommand{\arraystretch}{1.25}
\begin{array}{lll}
 D^2H(\varepsilon)(\sigma,\sigma)\geq h''(0)|\sigma|^2.
\end{array}
\eqno(2.6)
\]

In the proof of Theorem $1$, we need the following results.

The first results is Lemma $3.1$ in \cite{FZ}, which is
 a slight extension of a result of Giaquinta and Modica \cite{GM}. Define $Q_R(z_0)=\{(x,y)\in \mathbb{R}^2, \ |x-x_0|<R,  |y-y_0|<R\}$, where $,  z_0=(x_0, y_0)\in \mathbb{R}^2 $.

 \textbf{Lemma 2.1}\ \  Let $f,$ $f_1, \ldots,f_\ell$ denote non-negative functions from the space $L^1_{loc} (\mathbb{R}^2).$ Suppose
 further that we are given exponents $\alpha_1, \ldots, \alpha_\ell > 0.$ Then we can find a number $\delta_0 > 0$ depending on $\alpha_1, \ldots, \alpha_\ell$ 
as follows: if for $\delta \in (0, \delta_0)$ it is possible to calculate a constant $C (\delta) > 0$ such that the inequality
\[
\int_{Q_R (z)} f dx \le \delta \int_{Q_{2R} (z)} f dx + C (\delta) \sum^\ell_{j = 1} R^{-\alpha_j} \int_{Q_{2R}(z)} f_j dx
\]
holds for any choice of $Q_R (z) \subset \mathbb{R}^2,$ then there is a constant $C$ with the property
\[
\int_{Q_R (z)} f dx \le C \sum^\ell_{j = 1} R^{-\alpha_j} \int_{Q_{2R}(z)} f_j dx
\]
for all squares $Q_R (z)$.

Next, we need  a standard result concerning the ``divergence equations'', see e.g. \cite{Ga1,Ga2} or \cite{La}.

 \textbf{Lemma 2.2}\ \  Consider a function $f \in L^2 (Q_R (z))$ such that $\int_{Q_R (z)} f dx = 0.$ Then there exists a field
$v \in W_{0}^{1,2}(Q_R (z), \mathbb{R}^2)$ and a constant $C$ independent of $Q_R (z)$ such that we have
$\text{div} \ v = f$ on $ Q_R (z)$ together with the estimate
\[
\int_{Q_R (z)} |Dv|^2 dx \le C \int_{Q_R (z)} f^2 dx \, .
\]

Finally,  we need the  $L^2$-variant of the classical  Korn inequality.

 \textbf{Lemma 2.3}\ \  There is a constant $C$ independent of $Q_R (z)$ such that for all  $v \in W_0^{1,2}(Q_R (z), \mathbb{R}^2)$ it holds
\[
\int_{Q_R (z)} |Dv|^2 dx \le C \int_{Q_R (z)} |\varepsilon(v)|^2 dx \, .
\]

\bigbreak
\begin{center}
\textbf{\Large {\S 3.\ The Proof Of Theorem 1}}
\end{center}
The following energy estimate was proved as  Lemma $4.1$ in \cite{Fu1}.

 \textbf{Proposition 3.1}\ \  Suppose that $u \in C^2(\mathbb{R}^2,\mathbb{R}^2)$ is a bounded solution to $(1.1)$, where $h$
 satisfies  $(A1)-(A3)$. Then it holds
\[
\renewcommand{\arraystretch}{1.25}
\begin{array}{lll}
\displaystyle \int_{B_t(x_0)}H(\varepsilon(u))dx \leq C(t+1)
\end{array}
\]
for all discs $B_t(x_0) \subset \mathbb{R}^2.$

In the above Proposition $3.1$, we can replace $u \in
C^2(\mathbb{R}^2,\mathbb{R}^2)$ by  $u \in
C^1(\mathbb{R}^2,\mathbb{R}^2)$ and replace the discs $B_t(x_0)$ by
the squares $Q_t(x_0)$.  Indeed, if $u\in
C^1(\mathbb{R}^2,\mathbb{R}^2),$  in view of the elliptic condition
$(2.5),$ the system is local uniformly elliptic, using the
difference quote technique, it follows that $u\in
W^{2,2}_{\text{loc}}(\mathbb{R}^2,\mathbb{R}^2)\bigcap
C^1(\mathbb{R}^2,\mathbb{R}^2)$,  see \cite{Fu2}.  Furthermore, it
follows that  $\pi\in
W^{1,2}_{\text{loc}}(\mathbb{R}^2,\mathbb{R})$. Consequently, as
mentioned in \cite{Fu1}, we can follow the proof of Lemma $4.1$ in
\cite{Fu1} to show that Lemma $3.1$ holds under the assumption $u
\in C^1(\mathbb{R}^2,\mathbb{R}^2)$.

Now we have  $u\in W^{2,2}_{\text{loc}}(\mathbb{R}^2,\mathbb{R}^2)\bigcap C^1(\mathbb{R}^2,\mathbb{R}^2)$
and $\pi\in W^{1,2}_{\text{loc}}(\mathbb{R}^2,\mathbb{R}).$  Under the boundedness condition of the solution,
we would  like to prove the validity of
\[
\renewcommand{\arraystretch}{1.25}
\begin{array}{lll}
\displaystyle\int_{\mathbb{R}^2} D^2H(\varepsilon(u))(\varepsilon(\partial_ku), \varepsilon(\partial_ku))dx < \infty.
\end{array}
\eqno(3.1)
\]

Note that from $(2.6)$ and $(3.1)$ it immediately follows that
\[
\renewcommand{\arraystretch}{1.25}
\begin{array}{lll}
\displaystyle\int_{\mathbb{R}^2}|D\varepsilon(u)|^2 dx < \infty.
\end{array}
\eqno(3.2)
\]

To prove $(3.1),$  we go back  to the system $(1.1).$ For
any $\eta\in C_0^{\infty}(\mathbb{R}^2)$ with  $0\leq \eta \leq 1,$
letting $\varphi_k=\partial_{k}u\eta^2,$ we multiply $(1.1)$ with
$\partial_k\varphi_k$ and use integration by parts to obtain
\[
\renewcommand{\arraystretch}{1.25}
\begin{array}{lll}
\displaystyle\int_{\mathbb{R}^2}\partial_k\sigma \colon \varepsilon(\varphi_k)dx-\displaystyle\int_{\mathbb{R}^2}D\pi \cdotp \partial_k\varphi_k dx
-\displaystyle\int_{\mathbb{R}^2}u^i\partial_iu\cdotp\partial_k\varphi_k dx =0,
\end{array}
\]
where $\sigma:=DH(\varepsilon(u)):=\frac{h'(|\varepsilon(u)|)}{|\varepsilon(u)|}\varepsilon(u).$

Recalling $\varphi=\partial_{k}u\eta^2$ and using integration by parts again, we obtain
\[
\renewcommand{\arraystretch}{1.25}
\begin{array}{lll}
\displaystyle\int_{\mathbb{R}^2}\partial_k\sigma \colon \varepsilon(\partial_{k}u)\eta^2dx&=\displaystyle\int_{\mathbb{R}^2}\sigma\colon \partial_{k}(D\eta^2\odot \partial_ku)dx+\displaystyle\int_{\mathbb{R}^2}\partial_k\pi\text{div}(\varphi_k)dx\\
&+\displaystyle\int_{\mathbb{R}^2}u^i\partial_iu\cdotp\partial_k\varphi_k dx=:I+II+III,
\end{array}
\eqno(3.3)
\]
where ${\odot}$ is the symmetric product of vectors.

In the following, we will deal with  all of the terms in  $(3.3).$ For the left hand side, we have
\[
\renewcommand{\arraystretch}{1.25}
\begin{array}{lll}
\displaystyle\int_{\mathbb{R}^2}\partial_k\sigma \colon \varepsilon(\partial_{k}u)\eta^2dx=\displaystyle\int_{\mathbb{R}^2} D^2H(\varepsilon(u))(\varepsilon(\partial_ku), \varepsilon(\partial_ku))\eta^2 dx.
\end{array}
\eqno(3.4)
\]

We will estimate the items $I$, $III$ in the same way as that in \cite{Fu1}. For the completeness, we include the proofs here.

For $I$,  using Young's inequality and the  estimates $(2.3)$ and $(2.5)$ we have for any $ \delta>0,$
\[
\renewcommand{\arraystretch}{1.25}
\begin{array}{lll}
I&=\displaystyle\int_{\mathbb{R}^2}\sigma\colon \partial_{k}(D(\eta^2)\odot \partial_ku)dx\\
 &\leq C \biggl\{\displaystyle\int_{\mathbb{R}^2}h'(|\varepsilon(u)|) |Du|(|D\eta|^2+|D^2\eta|) dx+\displaystyle\int_{\mathbb{R}^2}h'(|\varepsilon(u)|)|D\eta|\eta |D^2u|dx\biggr\} \\
&\leq \delta  \displaystyle\int_{\mathbb{R}^2}\frac{h'(|\varepsilon(u)|)}{|\varepsilon(u)|}|D\varepsilon(u)|^2\eta^2dx +C(\delta)\displaystyle\int_{\mathbb{R}^2}h'(|\varepsilon(u)|)|\varepsilon(u)||D\eta|^2 dx\\
&+C \displaystyle\int_{\mathbb{R}^2}|Du|^2(|D\eta|^2+|D^2\eta|)dx + C \displaystyle\int_{\mathbb{R}^2}h'(|\varepsilon(u)|)^2(|D\eta|^2+|D^2\eta|)\\
&\leq \delta  \displaystyle\int_{\mathbb{R}^2} D^2H(\varepsilon(u))(\partial_k\varepsilon(u), \partial_k\varepsilon(u))\eta^2 dx+C(\delta)\displaystyle\int_{\mathbb{R}^2}h(|\varepsilon(u)|)|D\eta|^2dx\\
&+ C \displaystyle\int_{\mathbb{R}^2}h'(|\varepsilon(u)|)^2(|D\eta|^2+|D^2\eta|)dx + C  \displaystyle\int_{\mathbb{R}^2}|Du|^2(|D\eta|^2+|D^2\eta|)dx,
\end{array}
\eqno(3.5)
\]
where the relation $|D^2u(x)|\leq C|D\varepsilon(u)(x)|$ is used in the second inequality.

 For $III$ we have identity
\[
\renewcommand{\arraystretch}{1.25}
\begin{array}{lll}
III&=\displaystyle\int_{\mathbb{R}^2}u^i\partial_iu^j\partial_k(\partial_ku^j\eta^2 )dx=-\displaystyle\int_{\mathbb{R}^2}\partial_k(u^i\partial_iu^j)\partial_ku^j\eta^2dx\\
&=-\displaystyle\int_{\mathbb{R}^2}\partial_ku^i\partial_iu^j\partial_ku^j\eta^2dx-\displaystyle\int_{\mathbb{R}^2}u^i\frac{1}{2}\partial_i(|\partial_ku^j|^2)\eta^2dx\\
&=\frac{1}{2}\displaystyle\int_{\mathbb{R}^2}|Du|^2 u\cdotp
D\eta^2dx,
\end{array}
\eqno(3.6)
\]
where we use the identity $\partial_ku^i\partial_iu^j\partial_ku^j=0$ for divergence free vector $u$ in $2D.$

Finally, we estimate $II$. Here we just use equation $(1.1)$ to replace $D\pi$.
\[
\renewcommand{\arraystretch}{1.25}
\begin{array}{lll}
II&=\displaystyle\int_{\mathbb{R}^2}\partial_k\pi\text{div}(\varphi_k)dx=\displaystyle\int_{\mathbb{R}^2}\partial_k\pi \partial_k u \cdotp D\eta^2dx\\
&=-\displaystyle\int_{\mathbb{R}^2}\sigma_{ik}\partial_i(\partial_ku\cdotp
D\eta^2)dx -\displaystyle\int_{\mathbb{R}^2}u^i\partial_iu^k
\partial_k u\cdotp D\eta^2dx.
\end{array}
\eqno(3.7)
\]

We estimate the first integral of $(3.7)$ in the same way as that in the above for $I$ we have
\[
\renewcommand{\arraystretch}{1.25}
\begin{array}{lll}
II& \leq \delta  \displaystyle\int_{\mathbb{R}^2} D^2H(\varepsilon(u))(\varepsilon(\partial_ku), \varepsilon(\partial_ku))\eta^2 dx\\
  &+C(\delta)\displaystyle\int_{\mathbb{R}^2}h(|\varepsilon(u)|)|D\eta|^2dx+\displaystyle\int_{\mathbb{R}^2}h'(|\varepsilon(u)|)^2(|D\eta|^2+|D^2\eta|)dx\\
&+\displaystyle\int_{\mathbb{R}^2}|Du|^2(|D\eta|^2+|D^2\eta|)dx+ C \displaystyle\int_{\mathbb{R}^2}|Du|^2|u||D\eta|dx.
\end{array}
\eqno(3.8)
\]

Combining $(3.3),$ $(3.4),$  $(3.5),$ $(3.6)$ and $(3.8)$ and choosing $\delta=\frac{1}{4}$ , we end up with
\[
\renewcommand{\arraystretch}{1.25}
\begin{array}{lll}
&\displaystyle\int_{\mathbb{R}^2} D^2H(\varepsilon(u))(\varepsilon(\partial_ku), \varepsilon(\partial_ku))\eta^2 dx\\
  &\leq C \biggl\{\displaystyle\int_{\mathbb{R}^2}h(|\varepsilon(u)|)|D\eta|^2dx
+\displaystyle\int_{\mathbb{R}^2}h'(|\varepsilon(u)|)^2(|D\eta|^2+|D^2\eta|)dx\\
&+\displaystyle\int_{\mathbb{R}^2}|Du|^2(|D\eta|^2+|D^2\eta|)dx
+ \displaystyle\int_{\mathbb{R}^2}|Du|^2|u||D\eta|dx\biggr \}.
\end{array}
\eqno(3.9)
\]

Now, choosing $\eta\in C_0^\infty(Q_{\frac{3}{2}R}(x_0))$ such that $\eta \equiv 1$ in $Q_R(x_0)$, $|D\eta|\leq \frac{4}{R}$, and $|D^2\eta|\leq \frac{16}{R^2}$  we obtain from $(3.9)$ that
\[
\renewcommand{\arraystretch}{1.25}
\begin{array}{lll}
\displaystyle\int_{Q_R(x_0)}wdx
 &\leq C \biggl\{ \displaystyle\frac{1}{R^2}\int_{Q_{\frac{3}{2}R}(x_0)}h(|\varepsilon(u)|)dx
+\displaystyle\frac{1}{R^2}\int_{Q_{\frac{3}{2}R}(x_0)}h'(|\varepsilon(u)|)^2dx\\
&+\displaystyle\frac{1}{R^2}\int_{Q_{\frac{3}{2}R}(x_0)}|Du|^2dx\biggr \}
+\displaystyle\frac{C(\|u\|_{L_\infty})}{R}\int_{T_{\frac{3}{2}R}(x_0)}|Du|^2dx,
\end{array}
\eqno(3.10)
\]
where $w:=D^2H(\varepsilon(u))(\varepsilon(\partial_ku), \varepsilon(\partial_ku))$ and $T_{\frac{3}{2}R}(x_0)=Q_{\frac{3}{2}R}(x_0)\backslash\overline{Q_R(x_0)}.$\\

We will show that it follows from $(3.10)$ that $(3.1)$ holds. The proof is the same as in \cite{Fu1}. For the completeness,
we conclude the proof here, some steps are a little bit different from those in \cite{Fu1}.

Let  $\xi\in C_0^\infty(Q_{2R}(x_0))$ be the cut-off function such that,  $0\leq \xi \leq 1, \ \xi\equiv 1\  \text{on}\  Q_{\frac{3}{2}R}(x_0)$ and $|D\xi|\leq \frac{4}{R}.$  We have by Lemma $2.3$ that
\[
\renewcommand{\arraystretch}{1.25}
\begin{array}{lll}
\displaystyle\int_{Q_{\frac{3}{2}R}(x_0)}|Du|^2dx &\leq \displaystyle\int_{Q_{2R}(x_0)}{\xi}^2|Du|^2dx \\
&\leq C \biggl(\displaystyle\int_{Q_{2R}(x_0)}|D(u\xi)|^2dx +\displaystyle\int_{Q_{2R}(x_0)}|u|^2|D\xi|^2dx\biggr )\\
&\leq C \biggr(\displaystyle\int_{Q_{2R}(x_0)}|\varepsilon(u\xi)|^2dx + \displaystyle\int_{Q_{2R}(x_0)}|u|^2|D\xi|^2dx\biggr)\\
&\leq C \biggl( \displaystyle\int_{Q_{2R}(x_0)}|\varepsilon(u)|^2dx+\displaystyle\frac{1}{R^2}\int_{Q_{2R}(x_0)}|u|^2dx \biggr ).
\end{array}
\eqno(3.11)
\]

Now, since $\frac{h'(t)}{t}$ is an increasing function and $h$ satisfies $(2.3),$ for any $L>0,$ it follows that
\[
\renewcommand{\arraystretch}{1.25}
\begin{array}{lll}
\displaystyle\int_{Q_{\frac{3}{2}R}(x_0)}h'(|\varepsilon(u)|)^2 dx \\
\leq \displaystyle\int\limits_{\{x\in Q_{\frac{3}{2}R}(x_0),\  |\varepsilon(u)| \leq L\}}h'(L)^2 dx + \displaystyle\int\limits_{\{x\in Q_{\frac{3}{2}R}(x_0),\  |\varepsilon(u)|> L\}}h'(|\varepsilon(u)|)^2dx\\
\leq  \displaystyle C R^2 h'(L)^2+ \displaystyle \frac{C}{L^2} \int_{Q_{\frac{3}{2}R}(x_0)}h(|\varepsilon(u)|)^2 dx.
\end{array}
\eqno(3.12)
\]

By Sobolev inequality  we have
\[
\renewcommand{\arraystretch}{1.25}
\begin{array}{lll}
\displaystyle\int_{Q_{\frac{3}{2}R}(x_0)}h(|\varepsilon(u)|)^2 dx \leq \displaystyle C\biggl(\int_{Q_{2R}(x_0)}|D(\xi h(|\varepsilon(u)|))| dx\biggr)^2\\
\leq C\biggl(\displaystyle  \int_{Q_{2R}(x_0)}|D\xi| h(|\varepsilon(u)|)dx\biggr)^2+ C\biggl(\displaystyle  \int_{Q_{2R}(x_0)}\xi h'(|\varepsilon(u)|)|D\varepsilon(u)| dx\biggr)^2.
\end{array}
\eqno(3.13)
\]

To estimate the second term in the above inequality, by H\"{o}lder inequality we have
\[
\renewcommand{\arraystretch}{1.25}
\begin{array}{lll}
\biggl(\displaystyle  \int_{Q_{2R}(x_0)}\xi h'(|\varepsilon(u)|)|D\varepsilon(u)| dx\biggr)^2\\
 \displaystyle \leq \int_{Q_{2R}(x_0)}\xi h'(|\varepsilon(u)|)|\varepsilon(u)|dx\displaystyle\int_{Q_{2R}(x_0)}\xi\frac{h'(|\varepsilon(u)|)}{|\varepsilon(u)|}|D\varepsilon(u)|^2dx.
\end{array}
\eqno(3.14)
\]

Combining  the estimates $(3.13)$ and $(3.14)$ and recalling  the definitions of $\xi$ and $w$ we have

\[
\renewcommand{\arraystretch}{1.25}
\begin{array}{lll}
&\displaystyle\int_{Q_{\frac{3}{2}R}(x_0)}h(|\varepsilon(u)|)^2 dx \\
&\leq \displaystyle C \frac{1}{R^2}\biggl(\displaystyle  \int_{Q_{2R}(x_0)} h(|\varepsilon(u)|)dx\biggr)^2+\displaystyle C\int_{Q_{2R}(x_0)} h(|\varepsilon(u)|)dx\int_{Q_{2R}(x_0)} w dx.
\end{array}
\eqno(3.15)
\]

Thus  $(3.12)$ and $(3.15)$  give us
\[
\renewcommand{\arraystretch}{1.25}
\begin{array}{lll}
\displaystyle\int_{Q_{\frac{3}{2}R}(x_0)}h'(|\varepsilon(u)|)^2 dx&\leq \displaystyle C R^2 h'(L)^2+C\frac{1}{L^2}\frac{1}{R^2}  \biggl(\int_{Q_{2R}(x_0)}h(|\varepsilon(u)|) dx\biggr)^2\\
&+\displaystyle C \frac{1}{L^2} \int_{Q_{2R}(x_0)}h(|\varepsilon(u)|) dx \int_{Q_{2R}(x_0)}w dx.
\end{array}
\eqno(3.16)
\]

Then it follows from $(3.10), (3.11)$ and $(3.16)$ that
\[
\renewcommand{\arraystretch}{1.25}
\begin{array}{lll}
\displaystyle\int_{Q_{R}(x_0)}wdx &\displaystyle\leq Ch'(L)^2 +\displaystyle C \frac{1}{L^2}\frac{1}{R^2}\int_{Q_{2R}(x_0)}h(|\varepsilon(u)|) dx  \int_{Q_{2R}(x_0)}w dx\\
&+\displaystyle C \frac{1}{L^2}\frac{1}{R^4} \biggl(\int_{Q_{2R}(x_0)}h(|\varepsilon(u)|) dx\biggr)^2+\displaystyle C\frac{1}{R^2}  \int_{Q_{2R}(x_0)}h(|\varepsilon(u)|) dx\\
&+\displaystyle C \frac{1}{R} \int_{Q_{2R}(x_0)}h(|\varepsilon(u)|) dx+\displaystyle C\frac{1}{R^4}\int_{Q_{2R}(x_0)}|u|^2 dx \\
&+ \displaystyle C\frac{1}{R^3}\int_{Q_{2R}(x_0)}|u|^2 dx.
\end{array}
\eqno(3.17)
\]

Now we choose $L=\frac{1}{R\beta}, \beta>0.$   Then $(3.17)$ gives
\[
\renewcommand{\arraystretch}{1.25}
\begin{array}{lll}
\displaystyle\int_{Q_{R}(x_0)}wdx &\displaystyle\leq Ch'(\frac{1}{R\beta})^2 +\displaystyle C \beta^2\int_{Q_{2R}(x_0)}h(|\varepsilon(u)|) dx  \int_{Q_{2R}(x_0)}w dx\\
&+\displaystyle C \frac{\beta^2}{R^2} \biggl(\int_{Q_{2R}(x_0)}h(|\varepsilon(u)|) dx\biggr)^2+\displaystyle C\frac{1}{R^2}  \int_{Q_{2R}(x_0)}h(|\varepsilon(u)|) dx\\
&+\displaystyle C \frac{1}{R} \int_{Q_{2R}(x_0)}h(|\varepsilon(u)|) dx+\displaystyle C\frac{1}{R^4}\int_{Q_{2R}(x_0)}|u|^2 dx \\
&+ \displaystyle C\frac{1}{R^3}\int_{Q_{2R}(x_0)}|u|^2 dx.
\end{array}
\eqno(3.18)
\]

Letting $\beta^2=\frac{\varepsilon}{R+1}, \varepsilon<1$ in $(3.18)$ and applying Proposition $3.1$ we obtain that
\[
\renewcommand{\arraystretch}{1.25}
\begin{array}{lll}
\displaystyle\int_{Q_{R}(x_0)}wdx &\displaystyle\leq \displaystyle C \varepsilon  \int_{Q_{2R}(x_0)}w dx+C h'\biggl(\frac{(R+1)^{\frac{1}{2}}}{{\varepsilon}^{\frac{1}{2}}R}\biggr)^2 \\
&+\displaystyle C\frac{1}{R^2} \int_{Q_{2R}(x_0)}h(|\varepsilon(u)|)
dx
+\displaystyle C \frac{1}{R} \int_{Q_{2R}(x_0)}h(|\varepsilon(u)|) dx\\
&+\displaystyle C\frac{1}{R^4}\int_{Q_{2R}(x_0)}|u|^2 dx+
\displaystyle C\frac{1}{R^3}\int_{Q_{2R}(x_0)}|u|^2 dx.
\end{array}
\eqno(3.19)
\]

Taking into account of $(2.4)$ and choosing $\varepsilon$ small enough such that $\delta:=C\varepsilon<\frac{1}{2}$ we have
\[
\renewcommand{\arraystretch}{1.25}
\begin{array}{lll}
\displaystyle\int_{Q_{R}(x_0)}wdx &\displaystyle\leq \displaystyle  \delta\int_{Q_{2R}(x_0)}w dx+\displaystyle C(1+\frac{1}{R^{2m}})\\
&+\displaystyle C\frac{1}{R^2} \int_{Q_{2R}(x_0)}h(|\varepsilon(u)|)
dx
+\displaystyle C \frac{1}{R} \int_{Q_{2R}(x_0)}h(|\varepsilon(u)|) dx\\
&+\displaystyle C\frac{1}{R^4}\int_{Q_{2R}(x_0)}|u|^2 dx+
\displaystyle C\frac{1}{R^3}\int_{Q_{2R}(x_0)}|u|^2 dx,
\end{array}
\eqno(3.20)
\]

which, by Lemma $2.1$, give us
\[
\renewcommand{\arraystretch}{1.25}
\begin{array}{lll}
\displaystyle\int_{Q_{R}(x_0)}wdx &\displaystyle\leq \displaystyle C(1+\frac{1}{R^{2m}})+\displaystyle C \frac{1}{R^2} \int_{Q_{2R}(x_0)}h(|\varepsilon(u)|) dx\\
&+\displaystyle C\frac{1}{R} \int_{Q_{2R}(x_0)}h(|\varepsilon(u)|)
dx
+\displaystyle C\frac{1}{R^4}\int_{Q_{2R}(x_0)}|u|^2 dx\\
& + \displaystyle C\frac{1}{R^3}\int_{Q_{2R}(x_0)}|u|^2 dx. \\
\end{array}
\eqno(3.21)
\]

We apply Proposition $3.1$ again and let $R\rightarrow \infty$ in
$(3.21)$. we conclude the proof of $(3.1)$.

In the rest of the proof  we show that
\[
\renewcommand{\arraystretch}{1.25}
\begin{array}{lll}
\displaystyle\int_{\mathbb{R}^2} w dx = 0.
\end{array}
\eqno(3.22)
\]

Then it follows from $(3.22)$ that  $D\varepsilon(u)=0.$  And hence $D^2u=0.$  so  $u$ is  affine. Because we assume $u\in L^{\infty}(\mathbb{R}^2,\mathbb{R}^2),$ the claim of Theorem $1$  follows.

Now it remains to prove $(3.22)$. Let
\[
\renewcommand{\arraystretch}{1.25}
\begin{array}{lll}
w_{\infty}:=\displaystyle\int_{\mathbb{R}^2} w dx.
\end{array}
\]

Let us return to inequality $(3.10).$  To prove $(3.22),$ we need a delicate estimate for
the integral $\int_{T_{\frac{3}{2}R}(x_0)}|Du|^2 dx$. We choose a cut-off function $\zeta\in C_0^{\infty}(Q_{2R}(x_0)), 0\leq\zeta\leq 1, \zeta\equiv 1$ on $T_{\frac{3}{2}R}(x_0), \text{spt}\zeta \subset T_{2R}(x_0):=Q_{2R}(x_0)\backslash\overline{ Q_{\frac{R}{2}}(x_0)}$. By Lemma $2.3$,  we have
\[
\renewcommand{\arraystretch}{1.25}
\begin{array}{lll}
\displaystyle\int_{T_{\frac{3}{2}R}(x_0)}|Du|^2dx\leq \displaystyle\int_{Q_{2R}(x_0)}|D(u\zeta)|^2dx\leq C \displaystyle\int_{Q_{2R}(x_0)}|\varepsilon(u\zeta)|^2dx\\
\leq C \displaystyle\int_{Q_{2R}(x_0)}|u|^2|D\zeta|^2dx + C \displaystyle\int_{Q_{2R}(x_0)}\varepsilon_{ij}(u)\varepsilon_{ij}(u)\zeta^2dx.
\end{array}
\eqno(3.23)
\]

For the item $\int_{Q_{2R}(x_0)}\varepsilon_{ij}(u)\varepsilon_{ij}(u)\zeta^2dx$,  we
use the energy for the second order derivative to control it. By integration by parts and H\"{o}lder inequality we have
\[
\renewcommand{\arraystretch}{1.25}
\begin{array}{lll}
 &\displaystyle\int_{Q_{2R}(x_0)}\varepsilon_{ij}(u)\varepsilon_{ij}(u)\zeta^2dx\\
&=-\displaystyle\int_{Q_{2R}(x_0)} u^i\partial_j(\varepsilon_{ij}(u))\zeta^2dx-\displaystyle\int_{Q_{2R}(x_0)}u^i\varepsilon_{ij}(u)\partial_j\zeta^2dx\\
&\leq\biggl(\displaystyle\int_{Q_{2R}(x_0)}|u|^2dx\biggr)^{\frac{1}{2}}\biggl(\displaystyle\int_{Q_{2R}(x_0)}|D\varepsilon(u)|^2\zeta^4dx\biggr)^{\frac{1}{2}}\\
&+\biggl(\displaystyle\int_{Q_{2R}(x_0)}|u|^2|D\zeta^2|^2dx\biggr)^{\frac{1}{2}}\biggl(\displaystyle \int_{Q_{2R}(x_0)}|\varepsilon(u)|^2dx\biggr)^{\frac{1}{2}}.
\end{array}
\eqno(3.24)
\]

Putting together the estimates $(3.23)$ and $(3.24)$ and recalling the definition of $\zeta$ we obtain

\[
\renewcommand{\arraystretch}{1.25}
\begin{array}{lll}
\displaystyle\int_{T_{\frac{3}{2}R}(x_0)}|Du|^2dx
&\leq  C \displaystyle \frac{1}{R^2} \displaystyle\int_{Q_{2R}(x_0)}|u|^2dx \\
&+  \displaystyle C\frac{1}{R}\biggl(\displaystyle\int_{Q_{2R}(x_0)}|u|^2dx\biggr)^{\frac{1}{2}}\biggl(\displaystyle \int_{Q_{2R}(x_0)}|\varepsilon(u)|^2dx\biggr)^{\frac{1}{2}} \\
&+C\biggl(\displaystyle\int_{Q_{2R}(x_0)}|u|^2dx\biggr)^{\frac{1}{2}}\biggl(\displaystyle\int_{T_{2R}(x_0)}|D\varepsilon(u)|^2dx\biggr)^{\frac{1}{2}}.
\end{array}
\eqno(3.25)
\]

Then it follows from $(3.10)$, $(3.11)$, $(3.16)$ and $(3.25)$ that

\[
\renewcommand{\arraystretch}{1.25}
\begin{array}{lll}
\displaystyle\int_{Q_{R}(x_0)}wdx &\displaystyle\leq Ch'(L)^2 +\displaystyle C \frac{1}{L^2}\frac{1}{R^2}\int_{Q_{2R}(x_0)}h(|\varepsilon(u)|) dx  \int_{Q_{2R}(x_0)}w dx\\
&+\displaystyle C \frac{1}{L^2}\frac{1}{R^4} \biggl(\int_{Q_{2R}(x_0)}h(|\varepsilon(u)|) dx\biggr)^2+\displaystyle C\frac{1}{R^2}  \int_{Q_{2R}(x_0)}h(|\varepsilon(u)|) dx\\
&+\displaystyle C \frac{1}{R^2}\biggl(\int_{Q_{2R}(x_0)}|u|^2 dx\biggr)^\frac{1}{2} \biggl(\displaystyle\int_{Q_{2R}(x_0)}h(|\varepsilon(u)|) dx\biggr)^\frac{1}{2}\\
&+\displaystyle C\frac{1}{R}\biggl(\int_{Q_{2R}(x_0)}|u|^2 dx\biggr)^{\frac{1}{2}} \biggl(\displaystyle\int_{T_{2R}(x_0)}|D\varepsilon(u)|^2dx\biggr)^{\frac{1}{2}}\\
&+ \displaystyle C\frac{1}{R^3}\int_{Q_{2R}(x_0)}|u|^2 dx+\displaystyle C\frac{1}{R^4}\int_{Q_{2R}(x_0)}|u|^2 dx,
\end{array}
\eqno(3.26)
\]
from which, repeating the above steps $(3.18)$, $(3.19)$ and $(3.20)$ it follows that

\[
\renewcommand{\arraystretch}{1.25}
\begin{array}{lll}
\displaystyle\int_{Q_{R}(x_0)}wdx &\leq \displaystyle \frac{1}{2}\int_{Q_{2R}(x_0)}wdx+\displaystyle Ch'\biggl(\frac{(R+1)^{\frac{1}{2}}}{{\varepsilon}^{\frac{1}{2}}R}\biggr)^2 +\displaystyle C\frac{1}{R^2}  \int_{Q_{2R}(x_0)}h(|\varepsilon(u)|) dx\\
&+\displaystyle C \frac{1}{R^2}\biggl(\int_{Q_{2R}(x_0)}|u|^2 dx\biggr)^\frac{1}{2} \biggl(\displaystyle\int_{Q_{2R}(x_0)}h(|\varepsilon(u)|) dx\biggr)^\frac{1}{2}\\
&+\displaystyle C\frac{1}{R}\biggl(\int_{Q_{2R}(x_0)}|u|^2 dx\biggr)^{\frac{1}{2}} \biggl(\displaystyle\int_{T_{2R}(x_0)}|D\varepsilon(u)|^2dx\biggr)^{\frac{1}{2}}\\
&+ \displaystyle C\frac{1}{R^3}\int_{Q_{2R}(x_0)}|u|^2 dx+\displaystyle C\frac{1}{R^4}\int_{Q_{2R}(x_0)}|u|^2 dx,
\end{array}
\eqno(3.27)
\]
now let $R$ go to infinity. Then from which, together with the
boundedness  condition of  the solution $u$,   $h'(0)=0$, $(3.2)$
and Proposition $3.1$, gives us
\[
\renewcommand{\arraystretch}{1.25}
\begin{array}{lll}
w_{\infty}\leq \frac{1}{2} w_{\infty}.
\end{array}
\]

Thus $w_{\infty}=0$.   The proof is complete.

 \textbf{Acknowledgement}\ \ The author wants to thank Martin Fuchs for many discussions. This author was supported by the Academy of Finland.

\bibliographystyle{alpha}
\bibliography{Liouville1}

\end{document}